\documentclass{amsart}

\usepackage{amsmath, amssymb, mathrsfs, dsfont}

\numberwithin{equation}{section}

\begin{document}

\title[Uniqueness of Rankin-Selberg products]{Uniqueness of Rankin-Selberg products}


\author{Guy Henniart}
\address{Guy henniart \\ Institut Universitaire de France et Univ. Paris-Sud \\ Laboratoire de Math\'ematiques d'Orsay \\ CNRS \\ Orsay cedex F-91405, France}
\email{Guy.Henniart@math.u-psud.fr}

\author{Luis Lomel\'i}
\address{Luis Lomel\'i \\ University of Oklahoma \\ Department of Mathematics \\ Norman, OK 73019-3103}
\email{lomeli@math.ou.edu}


\date{February 2013}

\begin{abstract}
In the present paper, we show the equality of the $\gamma$-factors defined by Jacquet, Piatetski-Shapiro and Shalika with those obtained via the Langlands-Shahidi method. Contrary to the local proof given by Shahidi, our proof uses a refined version of the local-global principle for ${\rm GL}_n$ in positive characteristic, which has independent interest. The comparison of $\gamma$-factors is made via a uniqueness result for Rankin-Selberg $\gamma$-factors over a non-Archimedean local field of positive characteristic.
\end{abstract}

\maketitle

\section{introduction}

Let $F$ be a non-Archimedean locally compact field. The local Langlands conjecture for ${\rm GL}_n$ \cite{ht',h'00,lrs'} is known to preserve $L$-functions and $\varepsilon$-factors for pairs. Indeed, the family of correspondences when $n$ varies is characterized by such a preservation \cite{h'93}. In two previous papers the authors showed that the higher $L$-functions and $\varepsilon$-factors corresponding to the symmetric square, exterior square and Asai representations are preserved for $F$ of positive characteristic \cite{hl'11,hl'12}. Those factors are defined via the Langlands-Shahidi method in \cite{lomelipreprint,sha'90}.

Our proofs in \cite{hl'11,hl'12} are local-global and use the global Langlands correspondence proved by L. Lafforgue \cite{laff'}. In fact, we give a characterization of $\gamma$-factors by a series of local properties, combined with their occurrence in the global functional equation. The higher $\gamma$-factors in \cite{hl'11,hl'12} and the ones arising through the local Langlands correspondence \cite{lrs'} both satisfy these properties, hence are equal.

One of the local properties is the important multiplicativity property, which expresses the behavior of $\gamma$-factors under parabolic induction. It is in this property that in \cite{hl'11,hl'12} the local factors for pairs mentioned above make an appearance, but as defined by the Langlands-Shahidi method. On the other hand, Lafforgue uses the Rankin-Selberg factors defined by Jacquet, Piatetski-Shapiro and Shalika \cite{jpss'}. In \cite{sha'84}, Shahidi gave a proof that the two definitions give the same result; his proof is local in nature. In the present paper, we prove that our local-global approach gives a rather easy proof of Shahidi's result when $F$ has positive characteristic. Moreover, here we do not use Lafforgue's results --which would be unnatural as both types of factors are defined using only representations of linear groups ${\rm GL}_n(F)$ and not the Galois side of the Langlands correspondence. Instead, using only methods of Representation Theory of $\mathfrak{p}$-adic Reductive Groups and Automorphic Forms, we globalize a given local cuspidal representation of ${\rm GL}_n$ in an automorphic representation with controlled ramification at other places (Theorem~\ref{localglobalthm}). A variant of our result includes an automorphic analogue of the result of Katz and Gabber \cite{katz'} that we use in \cite{hl'11}.

The authors would like to thank C. J. Bushnell, J. Cogdell, B. Gross, P. Kutzko, A. Roche and F. Shahidi for helpful mathematical communications. The second author would like to thank the Automorphic Forms and Representation Theory group at the University of Oklahoma for providing an interesting working environment while this article was written.

\section{A uniqueness theorem}

\subsection{Notation} Let us adapt the notation of \cite{hl'11} and \S~6.1 of \cite{lomelipreprint} in order to better deal with $\gamma$-factors for pairs of representations of ${\rm GL}_m$ and ${\rm GL}_n$. Let $\mathscr{L}$ be the class of quadruples $(F,\pi_1,\pi_2,\psi)$ consisting of: a non-Archimedean local field $F$; smooth irreducible representations $\pi_1$ of ${\rm GL}_m(F)$ and $\pi_2$ of ${\rm GL}_n(F)$ (call $(m,n)$ the degree of the quadruple); and, a non-trivial character $\psi$ of $F$. We call a quadruple $(F,\pi_1,\pi_2,\psi) \in \mathscr{L}$ generic (resp. tempered, cuspidal) if the representations $\pi_1$ and $\pi_2$ are generic (resp. tempered, cuspidal). We fix a prime number $p$ and let $\mathscr{L}(p)$ be the class consisting of all $(F,\pi_1,\pi_2,\psi) \in \mathscr{L}$, with $F$ of characteristic $p$. Note that our proof is for $\mathscr{L}(p)$; the case of characteristic zero is mentioned in the remark following the proof of the theorem in \S~3.
 
Given a local non-Archimedean field $F$, we let $\mathcal{O}_F$ denote its ring of integers, $\mathfrak{p}_F$ its maximal ideal, $q_F$ the cardinality of its residue field, and $\left| \cdot \right|_F$ its absolute value. Given a global field $K$ and a non-Archimedean valuation $v$ of $K$, we write $\mathcal{O}_v$ for the ring of integers of $K_v$; and similarly for $\mathfrak{p}_v$, $q_v$ and $\left| \cdot \right|_v$. The cardinality of the field of constants of a global function field $K$ is denoted by $q$. Given ${\bf G} = {\rm GL}_l$, we write ${\bf P} = {\bf M}{\bf N}$ for a parabolic subgroup consisting of upper triangular block matrices with Levi subgroup ${\bf M}$. We let ${\bf B} = {\bf T}{\bf U}$ be the Borel subgroup of upper triangular matrices with maximal torus ${\bf T}$ and unipotent radical ${\bf U}$. Given a representation $\rho$, we let $\widetilde{\rho}$ denote its contragredient representation.

\subsection{Theorem.}\label{uniquenessthm} \emph{A rule $\gamma$ which assigns a rational function $\gamma(s,\pi_1 \times \pi_2,\psi) \in \mathbb{C}(q_F^{-s})$ to every $(F,\pi_1,\pi_2,\psi) \in \mathscr{L}(p)$, is uniquely determined by the following properties:}
\begin{itemize}
   \item[(i)] (Naturality) \emph{Let $(F,\pi_1,\pi_2,\psi) \in \mathscr{L}(p)$, and let $\eta : F' \rightarrow F$ be an isomorphism of local fields. Let $(F',\pi_1',\pi_2',\psi') \in \mathscr{L}(p)$ be the quadruple obtained from $(F,\pi_1,\pi_2,\psi)$ via $\eta$. Then
   \begin{equation*}
      \gamma(s,\pi_1 \times \pi_2,\psi) = \gamma(s,\pi_1' \times \pi_2',\psi').
   \end{equation*}
   }
  \item[(ii)] (Isomorphism). \emph{Let $(F,\pi_1,\pi_2,\psi) \in \mathscr{L}(p)$. If $(F,\pi_1',\pi_2',\psi) \in \mathscr{L}(p)$ is such that $\pi_i' \simeq \pi_i$, for $i=1$, $2$, then
     \begin{equation*}
        \gamma(s,\pi_1' \times \pi_2',\psi) = \gamma(s,\pi_1 \times \pi_2,\psi).
     \end{equation*}
  }
  \item[(iii)] (Compatibility with Tate's thesis). \emph{Let $(F,\chi_1,\chi_2,\psi) \in \mathscr{L}(p)$ be of degree $(1,1)$. Then
     \begin{equation*}
        \gamma(s,\chi_1 \times \chi_2,\psi) = \gamma(s,\chi_1 \chi_2,\psi),
     \end{equation*}
     where the right hand side is defined in Tate's thesis \cite{tate'}.
  }
  \item[(iv)] (Dependence on $\psi$). \emph{Let $(F,\pi_1,\pi_2,\psi) \in \mathscr{L}(p)$ be of degree $(m,n)$. Given $a \in F^\times$, let $\psi^a$ be the character of $F$ defined by $\psi^a(x) = \psi(ax)$. Then
     \begin{equation*}
        \gamma(s,\pi_1 \times \pi_2,\psi^a) = \omega_{\pi_1}(a)^{m} \omega_{\pi_2}(a)^{n} \left| a \right|_F^{mn(s-\frac{1}{2})} \gamma(s,\pi_1 \times \pi_2,\psi).
     \end{equation*}
  }
  \item[(v)] (Multiplicativity). \emph{For $1 \leq i \leq d_1$ and $1 \leq j \leq d_2$, let $(F,\tau_{1i},\tau_{2j},\psi) \in \mathscr{L}(p)$. For $h=1$,$2$, let $\pi_h$ be an irreducible subquotient of the representation of ${\rm GL}_{n_h}(F)$ parabolically induced from $\tau_{h1} \otimes \cdots \otimes \tau_{hd_h}$. Assume that for each $h =1$, $2$, either:
  \begin{itemize}
     \item[(a)] $\pi_h$ is generic or 
     \item[(b)] all of the $\tau_{hi}$'s are quasi-tempered and $\pi_h$ is the Langlands quotient of the parabolically induced representation.
  \end{itemize}
Then
     \begin{equation*}
        \gamma(s,\pi_1 \times \pi_2,\psi) = \prod_{i,j} \gamma(s,\tau_{1i} \times \tau_{2j},\psi).
     \end{equation*}
  }
  \item[(vi)] (Twists by unramified characters). \emph{Let $(F,\pi_1,\pi_2,\psi) \in \mathscr{L}(p)$, then
     \begin{equation*}
        \gamma(s+s_0,\pi_1 \times \pi_2,\psi) = \gamma(s, \left| \det(\cdot) \right|_F^{s_0} \pi_1 \times \pi_2,\psi).
     \end{equation*}
  }
  \item[(vii)] (Global functional equation). \emph{Let $K$ be a global function field of characteristic $p$. Let $\Psi = \otimes_v \Psi_v$ be a non-trivial character of $K \backslash \mathbb{A}_K$. Given cuspidal automorphic representations $\Pi_1 = \otimes' \, \Pi_{1,v}$ of ${\rm GL}_{n_1}(\mathbb{A}_K)$ and $\Pi_2 = \otimes' \, \Pi_{2,v}$ of ${\rm GL}_{n_2}(\mathbb{A}_K)$, let $S$ be a finite set of places of $K$ such that $\Pi_{1,v}$, $\Pi_{2,v}$ and $\Psi_v$ are unramified for $v \notin S$. Then
     \begin{equation*}
        L^S(s,\Pi_1 \times \Pi_2) = \prod_{v \in S} \gamma(s,\Pi_{1,v} \times \Pi_{2,v},\Psi_v) L^S(1-s,\widetilde{\Pi}_1 \times \widetilde{\Pi}_2).
     \end{equation*}
  }
\end{itemize} 

We provide the proof in \S~3. Before that, let us make a few remarks and derive a couple of consequences.

\subsection{Remark} In the global functional equation, partial $L$-functions are a product of local factors
\begin{equation*}
   L^S(s,\Pi_1 \times \Pi_2) = \prod_{v \notin S} L(s,\Pi_{1,v} \times \Pi_{2,v}).
\end{equation*}
More precisely, this product converges for $\Re(s)$ large enough. The resulting $L$-function has a meromorphic continuation to the complex $s$-plane and is a rational function on $q^{-s}$. The functions $L^S(s,\Pi_1 \times \Pi_2)$ and $L^S(s,\widetilde{\Pi}_1 \times \widetilde{\Pi}_2)$ verify property (vii).

\subsection{Remark}\label{multremark} Property (v) readily implies a stronger multiplicativity property: using the Langlands-Zelevinsky classification \cite{z'} we deduce the fact that if all the $\tau_{hi}$'s are cuspidal, then multiplicativity holds for all choices of subquotients $\pi_1$, $\pi_2$. In other words, the $\gamma$-factor $\gamma(s,\pi_1 \times \pi_2,\psi)$ depends only on the cuspidal supports of $\pi_1$, $\pi_2$, and the multiplicative property expresses $\gamma(s,\pi_1 \times \pi_2,\psi)$ as a product of $\gamma$-factors corresponding to those cuspidal supports. We then conclude that, with no special condition on the $\tau_{hi}$'s, multiplicativity holds for any choice of subquotients $\pi_1$, $\pi_2$. If we think about the Langlands correspondence, this corresponds to the fact that $\gamma$-factors for representations of the Weil-Deligne group of $F$ only depend on the underlying Weil group representation and are multiplicative with respect to direct sums.

\subsection{}\label{localfe} A first consequence is that $\gamma$-factors satisfy a local functional equation.

\medskip

\noindent{\bf Corollary.} A rule $\gamma$ as in Theorem~\ref{uniquenessthm} also satisfies the following property:
\begin{itemize}
   \item[(viii)] (Local functional equation). Let $(F,\pi_1,\pi_2,\psi) \in \mathscr{L}(p)$, then
   \begin{equation*}
      \gamma(s,\pi_1 \times \pi_2,\psi) \gamma(s,\widetilde{\pi}_1 \times \widetilde{\pi}_2,\overline{\psi}) = 1.
   \end{equation*}
\end{itemize}

\noindent{\emph{Proof}.} It is immediate that the rule $\gamma'$ on $\mathscr{L}(p)$ defined by $\gamma'(s,\pi_1 \times \pi_2,\psi) = \gamma(1-s,\widetilde{\pi}_1 \times \widetilde{\pi}_2,\overline{\psi})^{-1}$ satisfies properties (i) through (vii). \qed

\subsection{Generic $\gamma$-factors via different methods}\label{genericgamma} To derive our second consequence --the main reason for our investigation-- we consider two different rules on $\mathscr{L}(p)$. We begin by assuming that the underlying representations are generic, the general case is dealt with in \S~2.7.

First is the rule
\begin{equation*}
   (F,\pi_1,\pi_2,\psi) \mapsto \gamma(s,\pi_1 \times \pi_2,\psi)
\end{equation*}
defined by Jacquet, Piatetski-Shapiro and Shalika in \cite{jpss'}. Properties (i), (ii) and (vi) are easy consequences of the definitions. A direct proof of the dependance on $\psi$, Property~(iv), for Rankin-Selberg local factors is contained in the proof of Lemma~2.1 of \cite{cps'94}.
Multiplicativity for generic representations, Property~(v.a), is given by Theorem~3.1 of [\emph{op. cit.}]. The global functional equation --which we remark, involves only generic representations-- can be found as Theorem~2.3 of \cite{cps'}.

The second rule on $\mathscr{L}(p)$
\begin{equation*}
   (F,\pi_1,\pi_2,\psi) \mapsto \gamma(s,\pi_1 \otimes \widetilde{\pi}_2,\psi)
\end{equation*}
is obtained via the Langlands-Shahidi method \cite{lomelipreprint}. Let $(F,\pi_1,\pi_2,\psi) \in \mathscr{L}(p)$ be generic. Properties~(i), (ii), (vi) are immediate. Property~(iii) is given in Proposition~3.2 of [\emph{op. cit.}]. Whereas multiplicativity, Property~(v.a), can be found in equation~(6.5) of [\emph{op. cit.}]. We prove property~(iv) in the following lemma.

\medskip

\noindent{\bf Lemma.} \emph{Let $(F,\pi_1,\pi_2,\psi) \in \mathscr{L}(p)$ be generic of degree $(m,n)$. Given $a \in F^\times$, let $\psi^a$ be the character of $F$ defined by $\psi^a(x) = \psi(ax)$. Then
     \begin{equation*}
        \gamma(s,\pi_1 \otimes \widetilde{\pi}_2,\psi^a) = \omega_{\pi_1}(a)^{m} \omega_{\pi_2}(a)^{n} \left| a \right|_F^{mn(s-\frac{1}{2})} \gamma(s,\pi_1 \otimes \widetilde{\pi}_2,\psi).
     \end{equation*}
}
\emph{Proof.} The $\gamma$-factors are obtained via the Langlands-Shahidi method by considering ${\bf M} = {\rm GL}_m \times {\rm GL}_n$ as a maximal Levi subgroup of ${\bf G} = {\rm GL}_{m+n}$. Let $r$ denote the adjoint representation of ${}^LM$ on ${}^L\mathfrak{n}$. The character $\psi$ is used in \S~2.1 of \cite{lomelipreprint} to define a non-degenerate character of ${\bf U}(F)$, again denoted by $\psi$. We write $\psi_M$ for the restriction of $\psi$ to $U_M = {\bf M}(F) \cap {\bf U}(F)$. The characters $\psi$ and $\psi_M$ are then $w_0$-compatible in the notation of \S~6.2 of [\emph{op. cit.}]. We consider the $\psi_M$-generic representation $\pi = \pi_1 \otimes \widetilde{\pi}_2$ of ${\bf M}(F)$.

For $a \in F^\times$, let $t = {\rm diag}(a^{-(m+n-1)},a^{-(m+n-2)},\ldots,a,1)$. Let $\pi_t$ be the representation of $M$ given by $\pi_t(x) = \pi(t^{-1}xt)$. The character $\psi_t$ given by $\psi_t(u) = \psi(t^{-1}ut)$ is then obtained from $\psi^a$ and $\pi_t$ is $\psi^a_M$-generic. We can now explicitly apply (6.1) of [\emph{loc. cit.}] to the local coefficient in this setting:
\begin{align*}
   \gamma(s,\pi_1 \otimes \widetilde{\pi}_2,r,\psi^a) & = C'_{\psi_t}(s,\pi_t,w_0) \\
   									& = \omega_{{\pi}_1}(a)^{m} \omega_{\widetilde{\pi}_2}(a)^{-n} \left| a \right|_F^{mn(s-\frac{1}{2})} C'_\psi(s,\pi,w_0) \\
									& = \omega_{\pi_1}(a)^{m} \omega_{\pi_2}(a)^{n} \left| a \right|_F^{mn(s-\frac{1}{2})} \gamma(s,\pi_1 \otimes \widetilde{\pi}_2,r,\psi).
\end{align*}
\qed

Finally, we have a global functional equation: Let $K$ be a global function field of characteristic $p$, let $\Psi = \otimes_v \Psi_v$ be a non-trivial character of $k \backslash \mathbb{A}_K$, and let $\Pi_1$ and $\Pi_2$ be cuspidal automorphic representations of ${\rm GL}_{n_1}(\mathbb{A}_K)$ and ${\rm GL}_{n_2}(\mathbb{A}_K)$, respectively. Let $S$ be a finite set of places of $K$ such that $\Psi$ and $\Pi_i$, for $i=1$, $2$, are unramified outside of $S$. The global functional equation, Theorem~5.1 of [\emph{loc. cit.}] in its form of Property~6.5(vi) for $\gamma$-factors, gives
\begin{equation*}
   L^S(s,\Pi_1 \times \Pi_2) = \prod_{v \in S} \gamma(s,\Pi_{1,v} \otimes \widetilde{\Pi}_{2,v},\Psi_v) L^S(1-s,\widetilde{\Pi}_1 \times \widetilde{\Pi}_2).
\end{equation*}
Thus, Properties~(i) through (vii) hold for $\gamma(s,\pi_1 \otimes \widetilde{\pi}_2,\psi)$ and $\gamma(s,\pi_1 \times \pi_2,\psi)$, when $\pi_1$ and $\pi_2$ are generic.

\subsection{General case}\label{mainequality} When $\pi_1$ and $\pi_2$ are not necessarily generic, they can be written as Langlands' quotients. More specifically, $\pi_1$ and $\pi_2$ are quotients of the representations $\xi$ and $\tau$, respectively, given as follows:
\begin{equation*}
   \xi = {\rm ind}_{P'}^{{\rm GL}_m(F)}(\xi_1 \otimes \cdots \otimes \xi_d), \quad 
   \tau = {\rm ind}_{P''}^{{\rm GL}_n(F)}(\tau_1 \otimes \cdots \otimes \tau_e);
\end{equation*}
each $\xi_i = \left| \det(\cdot) \right|_F^{u_i} \xi_{i,0}$ and $\tau_{j} = \left| \det(\cdot) \right|_F^{v_j} \tau_{j,0}$ has $\xi_{i,0}$ and $\tau_{j,0}$ tempered. Then, the factors $\gamma(s,\pi_1 \otimes \widetilde{\pi}_2,\psi)$ are defined by equation~(7.5) of \cite{lomelipreprint}. That is, $\gamma$-factors are defined for all $(F,\pi_1,\pi_2,\psi) \in \mathscr{L}(p)$ by means of multiplicativity, Property~(v.b), and Property~(vi). They satisfy the equation
\begin{equation*}
   \gamma(s,\pi_1 \times \pi_2,\psi) = \prod_{i,j} \gamma(s+u_i+v_j,\xi_{i,0} \times \tau_{j,0},\psi).
\end{equation*}
Rankin-Selberg $\gamma$-factors for non-generic representations are defined in \cite{jpss'} via the exact same procedure.

It is easy to verify that properties (i) through (vii) of the theorem hold for both $\gamma(s,\pi_1 \otimes \widetilde{\pi}_2,r,\psi)$ and $\gamma(s,\pi_1 \times \pi_2,\psi)$. We have proved the following:

\subsection*{Corollary.} \emph{Let $(F,\pi_1,\pi_2,\psi) \in \mathscr{L}(p)$, then
\begin{equation*}
   \gamma(s,\pi_1 \times \pi_2,\psi) = \gamma(s,\pi_1 \otimes \widetilde{\pi}_2,\psi).
\end{equation*}
}

\subsection{Remark.} As mentioned in the introduction, this gives a local-global proof of a result due to Shahidi \cite{sha'84}; his proof is purely local. The idea of the local-global approach has a long history: one could trace it to the classical derivation of local class field theory from global class field theory. Already Deligne treats Artin $L$-functions and root numbers for Galois representations by incorporating twists by highly ramified characters. The corresponding stability property of the $\gamma$-factors $\gamma(s,(\pi_1 \otimes \eta) \times \pi_2,\psi)$, $\eta$ sufficiently ramified, can be found in \cite{js'85}. This property plays a dominant role in the proof of the local Langlands conjecture \cite{lrs',ht',h'00} and is also used in \cite{laff'}. In the Langlands-Shahidi method, it is a conjecture that the local coefficient is stable under highly ramified twists which has already been proved in several cases involving ${\rm GL}_n$ (further references may be found in \cite{sha'preprint}). Note, however, that in our proof in positive characteristic we do not use the stability property.

\subsection{$L$-functions and root numbers} We can define local $L$-functions and $\varepsilon$-factors via $\gamma$-factors via the relationship:
\begin{equation*}
   \gamma(s,\pi_1 \times \pi_2,\psi) = \varepsilon(s,\pi_1 \times \pi_2,\psi) \dfrac{L(1-s,\widetilde{\pi}_1 \times \widetilde{\pi}_2)}{L(s,\pi_1 \times \pi_2)}.
\end{equation*}
First for for tempered $(F,\pi_1,\pi_2,\psi) \in \mathscr{L}$ and then by means of Langlands classification and analytic continuation. Hence, we can use Corollary~\ref{mainequality} to obtain the following equality of local factors.

\medskip

\subsection*{Corollary.} \emph{Let $(F,\pi_1,\pi_2,\psi) \in \mathscr{L}(p)$. Then
\begin{align*}
   L(s,\pi_1 \times \pi_2) = L(s,\pi_1 \otimes \widetilde{\pi}_2) \ \text{ and } \ \varepsilon(s,\pi_1 \times \pi_2,\psi) = \varepsilon(s,\pi_1 \otimes \widetilde{\pi}_2,\psi).
\end{align*}
}

\section{Proof of Theorem~\ref{uniquenessthm}}

\subsection{} Assume that we have two rules $\gamma$ and $\gamma'$ on $\mathscr{L}(p)$ satisfying the conditions of Theorem~\ref{uniquenessthm}. We want to prove that 
\begin{equation*}
 \gamma(s,\pi_1 \times \pi_2,\psi) = \gamma'(s,\pi_1 \times \pi_2,\psi)
\end{equation*}
for all for all $(F,\pi_1,\pi_2,\psi) \in \mathscr{L}(p)$. We prove it by induction on $m + n$, where $(m,n)$ is the degree of $(F,\pi_1,\pi_2,\psi)$. The case of $m + n = 2$ is given by property~(iii). By Remark~\ref{multremark}, we have equality by induction if $\pi_1$ or $\pi_2$ is not cuspidal. So, we assume that $\pi_1$ and $\pi_2$ are cuspidal. By property~(vi), we may even assume that both are unitary. Note also that by property~(iv), it is enough to prove the equality for a fixed $\psi$ and equality follows for every $\psi$.

\subsection{} Let $k$ be the residue field of $F$, and $T$ be the usual choice of coordinate on the affine line $\mathbb{A}_k^1$ over $k$, so that $K = k(T)$ is the function field of the projective line $\mathbb{P}_k^1$ over $k$. By property~(i) we may assume that $F$ is the completion $K_0$ of $K$ at the point $0$, so that $\pi_1$ and $\pi_2$ are cuspidal unitary representations of ${\rm GL}_{m}(K_0)$ and ${\rm GL}_{n}(K_0)$, respectively. We choose a non-trivial additive character $\Psi$ of $\mathbb{A}_K / K$, and assume, as we may, that $\Psi_0 = \psi$.

Now we use the following local-global theorem, to be proved in \S~4 together with the variant mentioned in the introduction.

\subsection{Theorem}\label{localglobalthm} Let $\pi$ be a cuspidal unitary representation of ${\rm GL}_n(K_0)$. Then there exists a cuspidal unitary automorphic representation $\Pi = \otimes_v \Pi_v$ of ${\rm GL}_n(\mathbb{A}_K)$ whose local components $\Pi_v$ satisfy:
\begin{itemize}
   \item[(i)] $\Pi_0 \simeq \pi$;
   \item[(ii)] at places distinct from $0$, $1$ and $\infty$, $\Pi_v$ is unramified;
   \item[(iii)] $\Pi_1$ is a subquotient of an unramified principal series representation;
   \item[(iv)] $\Pi_\infty$ is a subquotient of a tamely ramified principal series representation.
\end{itemize}

\subsection{} Let $(F,\pi_1,\pi_2,\psi) \in \mathscr{L}(p)$ be cuspidal unitary. Applying the theorem to the representations $\pi_1$ and $\pi_2$, we obtain cuspidal automorphic representations $\Pi_1$ of ${\rm GL}_m(\mathbb{A}_K)$ and $\Pi_2$ of ${\rm GL}_n(\mathbb{A}_K)$.  Then, with the notation of Property~(vii), there is a finite set of places $S$ containing $0$ such that the global functional equation is satisfied by both $\gamma$ and $\gamma'$. Hence
\begin{equation*}
   \prod_{v \in S} \gamma(s,\Pi_{1,v} \times \Pi_{2,v},\Psi_v) = \prod_{v \in S} \gamma'(s,\Pi_{1,v} \times \Pi_{2,v},\Psi_v),
\end{equation*}
where $\Pi_{1,v}$ and $\Pi_{2,v}$ are principal series representations for $v \in S - \left\{ 0 \right\}$. Applying the already established non-cuspidal case at these places yields
\begin{equation*}
   \prod_{v \in S - \left\{ 0 \right\} } \gamma(s,\Pi_{1,v} \times \Pi_{2,v},\Psi_v) = \prod_{v \in S - \left\{ 0 \right\} } \gamma'(s,\Pi_{1,v} \times \Pi_{2,v},\Psi_v).
\end{equation*}
 The functions $L^S(s,\Pi_1 \times \Pi_2)$ and $L^S(s,\widetilde{\Pi}_1 \times \widetilde{\Pi}_2)$ appearing in the functional equation are uniquely determined. Hence, at the remaining place, we have
 \begin{equation*}
    \gamma(s,\Pi_{1,0} \times \Pi_{2,0},\Psi_0) = \gamma'(s,\Pi_{1,0} \times \Pi_{2,0},\psi).
\end{equation*}
The desired equality then follows from property~(ii).
 
\subsection{Remark.} Our proof of Corollary~\ref{mainequality} can be adapted to also work in characteristic $0$; we now provide a sketch. We use Proposition~5.1 of \cite{sha'90} as the link between the local and the global theory, and reduce the proof of Corollary~\ref{mainequality} to the case where $F$ is Archimedean. The theory for Archimedean local fields is studied in \cite{js',sha'85}. The complex case is easily reduced to ${\rm GL}_1$ by multiplicativity, and we use compatibility with Tate's thesis for ${\rm GL}_1$. When $F$ is real, we are reduced to ${\rm GL}_1$ or to the case where $\pi_1$ or $\pi_2$ are discrete series for ${\rm GL}_2(\mathbb{R})$. Even in the latter case, discrete series are components of principal series. We thus reduce $\gamma(s,\pi_1 \times \pi_2,\psi)$ as before to a product of $\Gamma$-factors for ${\rm GL}_1(\mathbb{R})$.

\section{Proof of Theorem~\ref{localglobalthm}}

We refine the argument of Appendix~1 of \cite{h'83}, which uses Poincar\'e series to construct the global cuspidal automorphic representation, now with controlled ramification in positive characteristic. It is interesting to note that the main idea was extended by Vign\'eras \cite{v'} to include generic representations of quasi-split classical groups; and, also for generic representations, Shahidi obtained a further refinement in the case of number fields \cite{sha'90}.

\subsection{} Given a place $v$ of the global field $K$, let $\mathcal{O}_v$ denote the ring of integers of $K_v$ and $\mathfrak{p}_v$ the maximal ideal. Set $G = {\rm GL}_n$ and let $Z$ denote its center. We set $\mathcal{K}_v = G(\mathcal{O}_v)$. We also write $G_v$ instead of ${\rm GL}_n(K_v)$.

We first construct a function $f = \otimes_v f_v$ on $G(\mathbb{A}_K)$. For $v \notin \left\{ 0, 1, \infty \right\}$, let $f_v$ be the characteristic function of $\mathcal{K}_v$. Let $f_1$ be the characteristic function of the Iwahori subgroup $\mathcal{I}_1$ of $G_1$ made out of matrices in $\mathcal{K}_1$ which are upper triangular modulo $\mathfrak{p}_1$. Similarly, we let $f_\infty$ be the characteristic function of the pro-$p$-Iwahori subgroup $\mathcal{I}_\infty^1$ of $G_\infty$ made out of matrices in $\mathcal{K}_\infty$ which are lower triangular unipotent modulo $\mathfrak{p}_\infty$.

The choice of $f_0$ is more involved. By chapter~6 of \cite{bk'}, there is a pair $(J,\lambda)$, where $J$ is a subgroup of $G_0$ which is open, contains the center $Z_0 = Z(K_0)$ of $G_0$, and is compact modulo $Z_0$; and, $\lambda$ is an irreducible smooth representation of $J$ such that $\pi$ is isomorphic to the representation obtained from $\lambda$ by compact induction from $J$ to $G_0$. By conjugation, we can assume that the maximal compact subgroup $J^0$ of $G_0$, made out of the elements whose determinant has absolute value $1$, is included in $\mathcal{K}_0$. As the central character $\omega_\pi$ of $\pi$ is assumed to be unitary, $\lambda$ is a unitarizable representation, and we can choose a non-zero coefficient $f_0$ of $\lambda$ verifying $f_0(g^{-1}) = \overline{f_0(g)}$ for $g \in J$. We extend $f_0$ by $0$ outside $J$, to get a function on $G_0$, still denoted $f_0$, and which is a coefficient of $\pi$.

\subsection{} The Poincar\'e series ${\rm P}f$ attached to $f$ is defined by
\begin{equation*}
   {\rm P}f(g) = \sum_{\gamma \in G(K)} f(\gamma g),
\end{equation*}
for $g \in G(\mathbb{A}_K)$. When $g$ lies in a compact subset of $G(\mathbb{A}_K)$ the sum is finite, hence ${\rm P}f$ is a continuous function on $G(\mathbb{A}_K)$. We can be more precise, let $\gamma$ belong to $G(K) \cap \left( J \times \prod_{v \neq 0} \mathcal{K}_v \right)$. Then, by the product formula, $\gamma$ at the place $v = 0$ belongs to $J^0$, hence to $\mathcal{K}_0$. It follows that $\gamma$ belongs to $G(k)$. Furthermore, we see that if $\gamma$ at the place $v = 1$ belongs to $\mathcal{I}_1$ and at the place $v = \infty$ to $\mathcal{I}_\infty^1$, then $\gamma = I_n$ is the identity matrix. Let
\begin{equation*}
   \mathfrak{A} = J \times \mathcal{I}_1 \times \mathcal{I}_\infty^1 \times \prod_{v \notin \left\{ 0,1,\infty \right\} }\mathcal{K}_v.
\end{equation*}
It follows that ${\rm P}f$ is the function on $G(K) \cdot \mathfrak{A}$, trivial on $G(K)$, and coinciding with $f$ on $\mathfrak{A}$.

We remark that
\begin{equation*}
   \mathbb{A}_K^\times = K^\times \cdot \left( K_0^\times \times (1 + \mathfrak{p}_\infty) \times \!\! \prod_{v \notin \left\{ 0,\infty \right\} } \mathcal{O}_v^\times \right).
\end{equation*}
So, there is a unique unitary character $\omega : K^\times \backslash \mathbb{A}_K^\times \rightarrow \mathbb{C}^\times$ which restricts to $\omega_\pi$ on $K_0^\times$; $\omega_{\pi_\infty}$ is trivial on $1+\mathfrak{p}_\infty$; and, $\omega_{\pi_v}$ is trivial on $\mathcal{O}_v^\times$ for $v \notin \left\{ 0, \infty \right\}$. By construction, ${\rm P}f$ transforms via $\omega$ under translation by $Z(\mathbb{A}_K) = \mathbb{A}_K^\times$.

\subsection{} Reasoning as in the Appendix, p.~147, of \cite{h'83}, $\left| {\rm P}f \right| \in L^2(Z(\mathbb{A}_K)G(K) \backslash G(\mathbb{A}_K))$; and, since $f_0$ is a coefficient of the cuspidal representation $\pi$, ${\rm P}f$ is a cuspidal function: it belongs to the space $L_0^2(G,\omega)$ of cuspidal automorphic functions on $G(K) \backslash G(\mathbb{A}_K)$ transforming via $\omega$ under the center.

The space $L_0^2(G,\omega)$ is an orthogonal sum of irreducible components, each occurring with multiplicity $1$. The projection onto any of those subspaces is $G(\mathbb{A}_K)$-equivariant. As ${\rm P}f$ is not $0$, we can choose such a component $\Pi$ such that the projection of ${\rm P}f$ on the space of $\Pi$ is not $0$ either.

Now $f$ transforms under the action of $\mathfrak{A} \cap Z(\mathbb{A}_k)$ via the restriction of $\omega$. Also, we obtain a function $\varphi$ of the Hecke algebra on $G(\mathbb{A}_K)$ by restriction of ${\rm P}f$ to $Z(\mathbb{A}_K) \cdot \mathfrak{A}$; then $\varphi$ transforms under $Z(\mathbb{A}_K)$ via $\omega$. It acts by convolution on $L_0^2(G,\omega)$. Then
\begin{align*}
   \left( \varphi * {\rm P}f \right) (1) & = \int_{Z(\mathbb{A}_K) \backslash G(\mathbb{A}_K)} \varphi(g^{-1}) {\rm P}f(g) \, dg \\
   				   & = \int_{\mathfrak{A} \cap Z(\mathbb{A}_K) \backslash \mathfrak{A}} \overline{f(g)} {\rm P}f(g) \, dg \neq 0,
\end{align*}
where the last integral being non-zero because ${\rm P}f = f$ on $\mathfrak{A}$. Hence, $\varphi$ does not annihilate ${\rm P}f$ and it follows that $\varphi$ does not annihilate $f$ either. Taking into account the definition of $f$ and $\varphi$, we see that: $\Pi_v$ has non-zero fixed vectors under $\mathcal{K}_v$, for $v \notin \left\{ 0, 1, \infty \right\}$; $\Pi_1$ has a non-zero fixed vector under $\mathcal{I}_1$; and $\Pi_\infty$ has a non-zero fixed vector under $\mathcal{I}_\infty^1$. Moreover, since $f_0$ is a coefficient of $\pi$, $\Pi_0$ is equivalent to $\pi$.

\subsection{} From \S~9.2 of \cite{bk'98}, the trivial character of the Iwahori subgroup $\mathcal{I}_1$ of $G_1$ is a type (in the sense of \S\S~3, 4 of [\emph{op. cit.}]) for the Berstein component of the trivial character of $G_1$; this implies that any smooth irreducible representation of $G_1$ with a non-zero fixed vector under $\mathcal{I}_1$ is a subquotient of an unramified principal series. This proves condition~(iii) of the theorem.

We now turn to condition~(iv) of the theorem. Let $\mathcal{I}_\infty$ denote the Iwahori subgroup consisting of matrices in $\mathcal{K}_\infty$ which are lower triangular ${\rm mod}~\mathfrak{p}_\infty$. If a smooth irreducible representation $\rho$ of $G_\infty$ has a non-zero fixed vector under $\mathcal{I}_\infty^1$, the space $V$ of such fixed vectors --which is finite-dimensional, $\rho$ being admissible-- is stable under $\mathcal{I}_\infty$. Consequently $V$ contains an irreducible representation of $\mathcal{I}_\infty$, which is trivial on $\mathcal{I}_\infty^1$. Since $\mathcal{I}_\infty / \mathcal{I}_\infty^1$ is abelian, such a representation has dimension~$1$. Hence, we conclude that $\rho$ contains a character $\chi$ of $\mathcal{I}_\infty$, trivial on $\mathcal{I}_\infty^1$.

By the work of L. Morris \cite{m'}, such a character $\chi$ of $\mathcal{I}_\infty$ is a type for $G_\infty$ --when $\chi$ is trivial, this gives the above mentioned result of \S~9.2 of \cite{bk'98}. More precisely, the compact open subgroup $\mathcal{I}_\infty$ is contained in $\mathbf{B}^-(K_\infty)$, where $\mathbf{B}^-$ is the Borel subgroup of lower triangular matrices with maximal torus $\mathbf{T}$ of diagonal matrices. The character $\chi$ of $\mathcal{I}_\infty$ restricts to a character $\chi_U$ of $\mathcal{I}_\infty \cap \mathbf{T}(K_\infty)$; and, \cite{m'} states that $(\mathcal{I}_\infty,\chi)$ is a cover for $(\mathcal{I}_\infty \cap \mathbf{T}(K_\infty),\chi_U)$ --in the sense of \S~8 of \cite{bk'98}. The character $\chi_U$ is a type for characters of ${\bf T}(K_\infty)$ restricting to $\chi_U$. It follows from \S~8.3 of [\emph{op. cit.}] that $\chi$ is a type for the corresponding Berstein component, that of the subquotients of principal series induced from such characters. Hence, condition~(iv) of the theorem is satisfied. \qed

\subsection{Remark.} A similar reasoning shows the existence of an automorphic cuspidal representation $\Pi$ of $G(\mathbb{A}_K)$ such that $\Pi_v$ is unramified for $v \notin \left\{ 0, \infty \right\}$, $\Pi_0$ is equivalent to $\pi$, and $\Pi_\infty$ has non-zero fixed points under $\mathcal{K}_\infty^1 = 1 + M_n(\mathfrak{p}_\infty)$. This provides the analogue of the result of Katz and Gabber used in \cite{hl'11,hl'12}.

\bigskip

\end{document}